\documentclass[11pt]{article}

\usepackage{amssymb}
\usepackage{amsmath}
\usepackage{amsthm}
\usepackage{graphicx}
\usepackage{accents}

\newcommand{\R}{\mathbb R}
\newcommand{\I}{\mathbb I}
\newcommand{\Ro}{\accentset{\circ}{\mathbb R}}
\newcommand{\bPi}{\bar{\Pi}}
\newcommand{\bY}{\bar{Y}}
\newcommand{\bU}{\bar{U}}
\newcommand{\bdelta}{\bar{\delta}}
\newcommand{\tU}{\tilde{U}}
\newcommand{\hU}{\hat{U}}
\newcommand{\cU}{\mathcal{U}}
\newcommand{\LX}{\mathcal{L}}
\newcommand{\cE}{\mathcal{E}}
\newcommand{\cR}{\mathcal{R}}

\newcommand{\bfL}{{\bf L}}
\newcommand{\bfC}{{\bf C}}
\newcommand{\bfBV}{{\bf BV}}
\newcommand{\ds}{\displaystyle}
\newcommand{\qedwhite}{\hfill \ensuremath{\Box}}

\newtheorem{remark}{Remark}[section]
\newtheorem{theorem}{Theorem}[section]
\newtheorem{lemma}{Lemma}[section]

\newtheorem{definition}{Definition}[section]

\usepackage{xcolor,todonotes}

\author{J.~Lang, P.~Mindt}
\title{Entropy-Preserving Coupling Conditions for One-dimensional Euler Systems at Junctions}
\author{
Jens Lang\footnote{corresponding author}
\hspace{0.2cm}and\hspace{0.2cm}Pascal Mindt\\ \\
{\small \it Technische Universit\"at Darmstadt} \\
{\small \it Dolivostra{\ss}e 15, 64293 Darmstadt, Germany}\\
{\small lang@mathematik.tu-darmstadt.de}\\
{\small mindt@mathematik.tu-darmstadt.de}
}
\date{September 28, 2017}
\begin{document}
\maketitle

\begin{abstract}
This paper is concerned with a set of novel coupling conditions
for the $3\times 3$ one-dimensional Euler system with source
terms at a junction of pipes with possibly different cross-sectional
areas. Beside
conservation of mass, we require the equality of the total enthalpy
at the junction and that the specific entropy for pipes
with outgoing flow equals the convex combination of all entropies
that belong to pipes with incoming flow.
Previously used coupling conditions include equality of pressure or
dynamic pressure. They are restricted to the special case of a junction
having only one pipe with outgoing flow direction. Recently,
{\sc Reigstad} [SIAM J. Appl. Math., 75:679--702, 2015] showed that
such pressure-based coupling conditions can produce non-physical
solutions for isothermal flows through the production of mechanical
energy. Our new coupling conditions ensure energy as well as entropy
conservation and also apply to junctions connecting an arbitrary number
of pipes with flexible flow directions. We prove the existence and
uniqueness of solutions to the generalised Riemann problem at a
junction in the neighbourhood
of constant stationary states which belong to the subsonic region. This
provides the basis for the well-posedness of the
homogeneous and inhomogeneous Cauchy problems for initial data with
sufficiently small total variation.
\end{abstract}



\section{Introduction}
We consider the one-dimensional polytropic Euler equations with source terms at a
network with one single junction connecting $N$ pipe sections of infinite length
\begin{eqnarray}
\label{euler_eqs}
\partial_t U^{(i)} + \partial_x F(U^{(i)}) &=& G(x,t,U^{(i)}),
\quad (x,t)\in\R^+\times\R^+,\\
\label{euler_init}
U^{(i)}(x,0) &=& U^{(i)}_0(x), \quad x\in \R^+,
\end{eqnarray}
for $i=1,\ldots,N$, with the thermodynamic variables and the flux functions
\begin{equation}
\label{euler_flux}
U^{(i)}=
\begin{pmatrix}
\rho_i\\[1mm]
\rho_iu_i\\[1mm]
E_i
\end{pmatrix}\quad\mbox{and}\quad
F(U^{(i)})=
\begin{pmatrix}
\rho_iu_i\\[1mm]
\rho_iu_i^2+p_i\\[1mm]
u_i(E_i+p_i)
\end{pmatrix}.
\end{equation}
Each pipe is described by a vector, $\nu_i\in\R^3\setminus\{0\}$,
originating from the common junction and parameterized by $x\in\R^+$,
the real halfline $[0,\infty)$. The
surface section of the pipe equals $\|\nu_i\|\!\ne\!0$. We assume
$\nu_i\ne \nu_j$ for $i\ne j$. Further, $\rho_i$ is
the density, $u_i$ is the velocity, $p_i$ is the pressure, and $E_i$ is
the total energy. The equation of state for an ideal polytropic gas in
the common form reads
\begin{equation}
E_i = \frac{p_i}{\gamma -1} + \frac12 \rho_iu_i^2
\end{equation}
with a suitable adiabatic exponent $\gamma>1$. For later use, we introduce
the mass flux, $q_i=\rho_iu_i$, the speed of sound, $c_i=\sqrt{\gamma p_i/\rho_i}$,
as well as the specific entropy $s_i$ and the total enthalpy $h_i$ defined by
\begin{equation}
s_i = c_v\,\ln\left( \frac{p_i}{\rho_i^\gamma}\right)\quad\mbox{and}\quad
h_i = \frac{E_i+p_i}{\rho_i}
\end{equation}
with the specific (constant) heat capacity $c_v>0$. More details about the underlying
thermodynamic principles can be found, e.g., in \cite[Sect.14.4]{LeVeque2002}. The
right-hand side vector $G(x,t,U^{(i)})$ describes source terms, e.g., gravity and friction.
We will first discuss the homogeneous case $G=0$, yielding a system of conservation laws
in (\ref{euler_eqs}), and extend our results to the inhomogeneous case later on through
operator splitting techniques, following known concepts.

The characteristic eigenvalues of the Euler equations are
\begin{equation}
\lambda_1(U) = u - c,\quad  \lambda_2(U) = u,\quad \lambda_3(U) = u + c.
\end{equation}
As usual in the literature, we also restrict our analysis to the subsonic region
defined by $|u|<c$, and introduce the two sets of subsonic data
\begin{eqnarray}
D^+ := \{ U=(\rho,\rho u,E)\in\Ro^+\times\R\times\Ro^+:
\,\lambda_1(U)<0<\lambda_2(U)<\lambda_3(U)\},\\
D^- := \{ U=(\rho,\rho u,E)\in\Ro^+\times\R\times\Ro^+:
\,\lambda_1(U)<\lambda_2(U)<0<\lambda_3(U)\},
\end{eqnarray}
with $\Ro^+=(0,\infty)$. Due to $\lambda_2(U)=u$ and the orientation of
the pipes, we can relate pipes with a flow direction towards the junction
with $D^-$ (incoming flow), while $D^+$ corresponds to pipes with flow
direction away from the junction (outgoing flow). The corresponding index
sets are defined by $\I_{i}:=\{i:U^{(i)}\in D^-\}$ and
$\I_{o}:=\{i:U^{(i)}\in D^+\}$. We will only consider cases with
$\I_i\cup\I_o=\{1,\ldots,N\}$.

The main challenge in network modelling is to prescribe a set of coupling
conditions  at the junction-pipe interfaces of the form
\begin{equation}
\label{euler_ccond}
\Phi \left( U^{(1)}(0^+,t),\ldots,U^{(N)}(0^+,t) \right) = \Pi(t),
\end{equation}
where $\Phi$ is a possibly nonlinear function of the traces
$U^{(i)}(0^+,t)=\lim_{x\rightarrow 0^+}U^{(i)}(x,t)$ of
the unknown variables and $\Pi$ is a coupling constant, which depends
only on time. The conditions are closely linked to the Euler equations
(\ref{euler_eqs}) and provide a relation between the flows in all pipes.
Various functions $\Phi$ have been proposed in the literature. We find
\[
\begin{array}{crl}
(M) & \sum_{i=1}^{N}\,\|\nu_i\|\,q_i(0^+,t) = 0,\;t>0
    & \mbox{(conservation of mass),}\\[2mm]
(E) & \sum_{i=1}^{N}\,\|\nu_i\|\,(u_i(E_i+p_i))(0^+,t) = 0,\;t>0
    & \mbox{(conservation of energy),}\\[2mm]
(P) & p_i(0^+,t) = p^*(t),\;t>0
    & \mbox{(equality of pressure),}\\[2mm]
(P_D)& (\rho_iu_i^2+p_i)(0^+,t) = P^*(t),\;t>0
    & \mbox{(equality of dynamic pressure),}\\[2mm]
(H)& h_i(0^+,t) = h^*(t),\;t>0
    & \mbox{(equality of enthalpy),}\\[2mm]
(S) & \sum_{i=1}^{N}\,\|\nu_i\|\,(q_is_i)(0^+,t) \ge 0,\;t>0
    & \mbox{(entropy increase),}\\[2mm]
\end{array}
\]
where $p^*(t)$, $P^*(t)$ and $h^*(t)$ are unique, scalar, momentum-
and enthalpy-related coupling constants, respectively. Note that
the dynamic pressure in $(P_D)$ equals the momentum flux in (\ref{euler_flux}).

{\sc Colombo} and {\sc Mauri} \cite{ColomboMauri2008} used coupling
conditions that include mass and energy conservation at the junction,
the equality of dynamic pressure as well as the entropy increase, i.e.,
the trace of the solution satisfies $(M)$, $(E)$, $(P_D)$, and $(S)$. They
proved the well-posedness of the Cauchy problem given by the equations
(\ref{euler_eqs}), (\ref{euler_init}), and (\ref{euler_ccond}) above,
under the standard condition that the total variation of the initial data
is sufficiently small. The proof was given for the special case of
$\I_o=\{1\}$ and $\I_i=\{2,\ldots,N\}$, i.e., one pipe with outgoing flow and
incoming flow in the remaining $N\!-\!1$ pipes. {\sc Herty} \cite{Herty2008}
replaced the coupling condition $(P_D)$ by the equality of pressure, $(P)$,
widely used in the engineering community to simulate gas networks. Following
the approach presented in \cite{ColomboMauri2008}, he also showed
well-posedness of the Cauchy problem for the special network studied there.
However, the comparison to two-dimensional
numerical results did not give a conclusion on whether dynamic pressure or
pressure is the most appropriate momentum-related coupling constant. The
one-dimensional coupling of two systems of Euler equations at a fixed
interface were studied by {\sc Chalons}, {\sc Raviart} and {\sc Seguin} in
\cite{ChalonsRaviartSeguin2008}. They discussed possible solutions to
coupled Riemann problems for three different types of coupling conditions.
{\sc Colombo} and {\sc Marcellini} \cite{ColomboMarcellini2010} investigated
the coupling of two pipes with different cross sectional areas and extended
their results to a more complex pipe with spatially varying cross sectional area.
An important and necessary assumption is the bounded total variation of
the pipe's area profile. Physically motivated coupling conditions for tunnel
fires in networks were formulated by {\sc Gasser} and {\sc Kraft}
\cite{GasserKraft2008}. They considered the small Mach number regime and
assumed a good mixing of the flow in the junction, which motivates
conservation of mass and internal energy, the equality of pressure and
an equal inflow condition for all densities of outgoing tunnels.

Pressure equality, $(P)$, as coupling condition for isothermal flow in pipeline
networks have been intensively studied by {\sc Banda}, {\sc Herty} and {\sc Klar}
\cite{BandaHertyKlar2006b,BandaHertyKlar2006a}.
Recently, {\sc Reigstad} \cite{Reigstad2015} (see also
\cite{MorinReigstad2015,Reigstad2014,ReigstadFlattenHaugenYtrehus2015})
showed for this type of flow that both coupling conditions $(P)$ and
$(P_D)$ deliver non-physical solutions
characterized by the production of mechanical energy at a junction
in a constructed test case with $N\!=\!3$. The main result of the
paper comprises the fact that only the Bernoulli invariant taken as
momentum-related coupling constant is proved to yield entropic solutions for
all subsonic flow conditions in the general case of a junction connecting $N$ pipes
of arbitrary cross-sectional area. The Bernoulli invariant equals the specific
stagnation enthalpy and thus can be seen as the enthalpy-related
coupling constants $h^*(t)$ in condition $(H)$ above. Together with
the conservation of mass and the relation $q_ih_i\!=\!u_i(E_i+p_i)$,
the equality of enthalpy at the junction immediately yields the conservation of
energy. Thus, $(M)$ and $(H)$ imply $(E)$ for the Euler system. In this sense,
the equality of enthalpy at the junction confirms the energy conservation there
and represents a first step towards answering the main question of how to close
the set of coupling conditions.

In contrast to the isothermal flow, the situation for the compressible Euler
equations with subsonic flow conditions is still unsettled and the analysis
suffers from the open question: What are further physically sound coupling
conditions for which well-posedness of Cauchy problems can be shown for the
general case of a junction connecting $N$ pipes of arbitrary cross-sectional
area and flexible flow directions? A common approach to tackle this question is
to consider a generalised Riemann problem at the junction. Suppose we ensure
mass conservation and the continuity of the enthalpy, i.e., $(M)$ and $(H)$
hold. Then, a closer inspection of the local solution structure of the
Riemann problem and the corresponding degrees of freedom (as done in
Sect.~\ref{sec:gen_riemann_prob})
shows that only one further coupling condition can be imposed for each of the
outgoing pipes. This observation also explains the choice of the special network in
\cite{ColomboMauri2008,Herty2008}. There, $(P)$ or $(P_D)$ were chosen instead
of $(H)$, and the conservation of energy was added, which allows to only consider one
outgoing pipe.

In this paper, we consider the equality of the entropy at the junction-pipe
interface for pipes with outgoing flow:
\begin{equation}
\begin{array}{crl}
(S_o) & s_i(0^+,t) = s^*(t),\;t>0,\;i\in\I_o
      & \mbox{(equality of outgoing entropy)},
\end{array}
\end{equation}
where the coupling constant $s^*(t)$ is identified as the convex combination of
all entropies that belong to the pipes with incoming flow. That is, we set
\begin{equation}
\begin{array}{crl}
\label{enthalpy_mix}
(S_i) & \ds s^*(t) = \frac{1}{\sum_{i\in\I_i}\,\|\nu_i\|q_i(0^+,t)}
        \,\sum_{i\in\I_i}\,\|\nu_i\|(q_is_i)(0^+,t)
      & \mbox{(entropy mix)}.
\end{array}
\end{equation}
Our choice is motivated by the assumption that gas flows entering a
junction mix perfectly, which was also used by {\sc Schmidt}, {\sc Steinbach},
and {\sc Willert} \cite{SchmidtSteinbachWillert2015} to derive a mixing temperature
at junctions and by {\sc Gasser} and {\sc Kraft} \cite{GasserKraft2008} to formulate
an equal inflow boundary condition for all densities of outgoing pipes. A
direct consequence of (\ref{enthalpy_mix}) and the conservation
of mass is the conservation of entropy per unit volume in smooth flows.
In this case, the
momentum equation in (\ref{euler_eqs}) can be equivalently reformulated to
$\partial_t(\rho s)+ \partial_x(qs)\!=\!0$ (see, e.g., \cite[Sect.14.5]{LeVeque2002}).
Thanks to (\ref{enthalpy_mix}), we have the identity
$\sum_{i\in\I_i}\,\|\nu_i\|(q_is_i)(0^+,t)=s^*(t)\sum_{i\in\I_i}\,\|\nu_i\|q_i(0^+,t)$,
and therefore
\begin{equation}
\begin{array}{crl}
(S') & \sum_{i=1}^{N}\,\|\nu_i\|\,(q_is_i)(0^+,t) = 0,\;t>0
     & \mbox{(entropy conservation}).
\end{array}
\end{equation}
The paper is organised as follows. In Sect.~\ref{sec:gen_riemann_prob},
we formulate the generalised Riemann problem at a junction with
the coupling conditions $(M)$, $(H)$, $(S_o)$, $(S_i)$ and show
its well-posedness. The corresponding Cauchy problem and its solution
are studied in Sect.~\ref{sec:cauchy_prob}. A summary is given in
Sect.~\ref{sec:summary}.

\section{Generalised Riemann problem at a junction}
\label{sec:gen_riemann_prob}
In this section, we show the well-posedness of the
coupling conditions $(M)$, $(H)$, $(S_o)$, and $(S_i)$ for
the homogeneous problem given by (\ref{euler_eqs}) with
$G\!=\!0$. To this end, we consider a generalised Riemann
problem at a junction and show that there exist a
unique self-similar solution in terms of the classical
Lax solution to standard Riemann problems. The theoretical
framework was introduced by {\sc Colombo} and {\sc Garavello}
\cite{ColomboGaravello2006} for the $p$-system and generalised
in \cite{ColomboMauri2008} to Euler systems.

Let denote by $\Omega_i=\{U^{(i)}\in\R^3:\rho_i>0,\,p_i>0\}$
nonempty sets and define the overall state space
$\Omega=\Omega_1\times\Omega_2\times\cdots\times\Omega_N$.
Furthermore, let $Y=(Y^{(1)},\ldots,Y^{(N)})$.

We first recall two basic definitions for generalised Riemann
problems at junctions.
\begin{definition}
The Riemann problem at a junction with $N$ pipes is defined through the set
of equations
\begin{equation}
\label{prob_jrm}
\begin{array}{rll}
\partial_t Y^{(i)} + \partial_x F(Y^{(i)}) &=& 0,
\quad (x,t)\in\R^+\times\R^+,\\[2mm]
\Phi\left( Y^{(1)}(0^+,t),\ldots,Y^{(N)}(0^+,t) \right) &=& \bPi,\\[2mm]
Y^{(i)}(x,0) &=& \bY^{(i)}_0, \quad x\in \R^+,
\end{array}
\end{equation}
for $i=1,\ldots,N$, where $\bY^{(1)}_0,\ldots,\bY^{(N)}_0$ are constant
thermodynamic states in $\Omega$ and $\bPi\in\R^d$ is also constant.
\end{definition}

\begin{definition}
\label{def_phi_sol}
A $\Phi$-solution to the Riemann problem (\ref{prob_jrm}) is a self-similar
function $Y(x,t):\R^+\times\R^+\rightarrow\Omega$ for which the following hold:
\begin{itemize}
\item[1.] There exists a constant state $Y_\ast(\bY_0)=\lim_{x\rightarrow 0^+}Y(x,t)$
such that all components $Y^{(i)}(x,t)$ coincide with the restriction to
$x>0$ of the Lax solution to the standard Riemann problem for $x\in\R$,
\begin{equation}
\label{prob_crm}
\begin{array}{rll}
\partial_t Y^{(i)} + \partial_x F(Y^{(i)}) &=& 0,\\[2mm]
Y^{(i)}(x,0) &=& \left\{
 \begin{array}{ll}
 \bY^{(i)}_0 & \mbox{if }x>0,\\[1mm]
 Y^{(i)}_\ast & \mbox{if }x<0.
 \end{array}
\right.
\end{array}
\end{equation}
\item[2.] The state $Y_\ast$ satisfies $\Phi(Y_\ast)=\bPi$
for all $t>0$.
\end{itemize}
\end{definition}
\begin{figure}[h]
\centering
\includegraphics[width=0.85\textwidth]{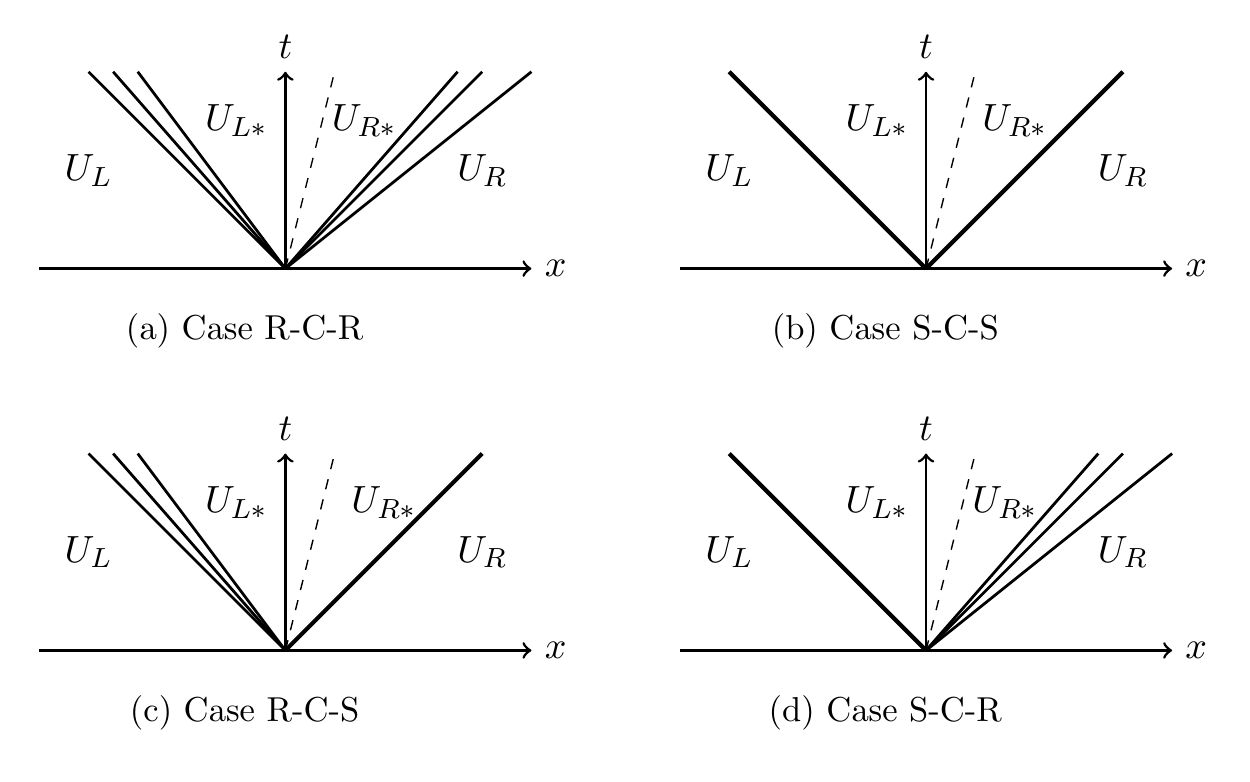}
\parbox{13cm}{
\caption{Possible wave patterns in the solution of Riemann problems
for the Euler equations: shock (S), contact (C) and rarefaction (R).}
\label{fig-4riemann}
}
\end{figure}
The solution of the standard Riemann problem (\ref{prob_crm})
with initial data $(U_L,U_R)$ for $x<0$ and $x>0$, respectively,
can be described by a set of elementary waves such as rarefaction, contact
and shock waves. The
three waves separate four constant states $(U_L,U_{L\ast},U_{R\ast},U_R)$.
The structure of the Euler equations reveals that the middle $2$-wave is
always a contact discontinuity while the left and right waves can be either
shock or rarefaction waves, see Fig.~\ref{fig-4riemann}. Further, both the
velocity and the pressure are constant across the contact discontinuity,
i.e., it holds
\begin{equation}
p_\ast = p_{L\ast}=p_{R\ast}\quad\mbox{ and }
\quad u_\ast = u_{L\ast}=u_{R\ast}.
\end{equation}
The four sought (constant) variables $(p_\ast,u_\ast,\rho_{L\ast},\rho_{R\ast})$
are implicitly defined by means of parametrisations of the Rankine-Hugoniot
jump condition and the Riemann invariants, see \cite[Sect.4]{Toro2009}
or \cite[Sect.14.11]{LeVeque2002} for more details. We have
\begin{eqnarray}
\label{eq_rm_u}
u_\ast = u_L-\psi(p_\ast,U_L) = u_R+\psi(p_\ast,U_R),\\
\label{eq_rm_rho}
\rho_{L\ast} = \phi(p_\ast,U_L),\quad \rho_{R\ast} = \phi(p_\ast,U_R),
\end{eqnarray}
where for $k=L,R$,
\begin{eqnarray}
\psi(p_\ast,U_k) &=&
\begin{cases}
  \ds \frac{2c_k}{\gamma-1}\left(\left(
    \frac{p_\ast}{p_k}\right)^{\frac{\gamma-1}{2\gamma}}-1\right)
    &\mbox{if }p_{\ast}\leq p_k \mbox{ (rarefaction)}\\[4mm]
  \ds (p_\ast-p_k)\,\left(\frac{1-\mu^2}{\rho_k(p_\ast+\mu^2p_k)}
    \right)^{\frac{1}{2}}
    &\mbox{if }p_{\ast}>p_k \mbox{ (shock)}
\end{cases}\\[2mm]
\phi(p_\ast,U_k) &=&
\begin{cases}
  \ds \rho_k\,\left(\frac{p_\ast}{p_k}\right)^{\frac{1}{\gamma}}
    &\mbox{if }p_{\ast}\leq p_k \mbox{ (rarefaction)}\\[4mm]
  \ds \rho_k\,\frac{p_\ast+\mu^2p_k}{\mu^2p_\ast+p_k}
    &\mbox{if }p_{\ast}>p_k \mbox{ (shock)}
\end{cases}
\end{eqnarray}
with $\mu^2=(\gamma-1)/(\gamma+1)$ and $c_k^2=\gamma p_k/\rho_k$.
Observe that the second equality in (\ref{eq_rm_u}) is used to
determine the parameter $p_\ast$. The functions $\psi(p_\ast,U_k)$
and $\phi(p_\ast,U_k)$ are twice continuously differentiable at
$p_\ast=p_k$. The total energy for the inner region can be
computed from $E_{k\ast}=p_\ast/(\gamma-1)+\rho_{k\ast}u_{\ast}^2/2$
for $k=L,R$.

For later use, let $\LX_1(\sigma,U_L)$ denote the $1$-Lax curve, which
parameterizes the $1$-wave curve through the state $U_L$
and describes all physical states on the right that can be reached from $U_L$
by either a shock wave for $\sigma>p_L$ or a rarefaction wave for
$\sigma\le p_L$. Using (\ref{eq_rm_u}) and (\ref{eq_rm_rho}), $\LX_1$ is
defined through
\begin{equation}
\label{eq_lax1}
\LX_1(\sigma,U_L) :=
\begin{pmatrix}
  \phi(\sigma,U_L)\\[1mm]
  \phi(\sigma,U_L)(u_L-\psi(\sigma,U_L))\\[1mm]
  \ds\frac{\sigma}{\gamma-1}+ \frac12\phi(\sigma,U_L)(u_L-\psi(\sigma,U_L))^2
\end{pmatrix}
\end{equation}
Analogously, let $\LX_3(\sigma,U_R)$ denote the $3$-Lax curve through
the state $U_R$, defined through
\begin{equation}
\LX_3(\sigma,U_R) :=
\begin{pmatrix}
  \phi(\sigma,U_R)\\[1mm]
  \phi(\sigma,U_R)(u_R+\psi(\sigma,U_R))\\[1mm]
  \ds\frac{\sigma}{\gamma-1}+ \frac12\phi(\sigma,U_R)(u_R+\psi(\sigma,U_R))^2
\end{pmatrix}
\end{equation}
We further recall the fact that for the $2$-contact discontinuity,
any state
\begin{equation}
\LX_2(\tau,\bar{U}) := \bar{U} + \tau
\left( 1,\,\bar{u},\,\frac{1}{2}\bar{u}^2\right)^T
\end{equation}
can be connected to $\bar{U}$ for sufficiently small $\tau\in\R$. This
defines the $2$-Lax curve.
\begin{figure}[h]
\centering
\includegraphics[width=0.9\textwidth]{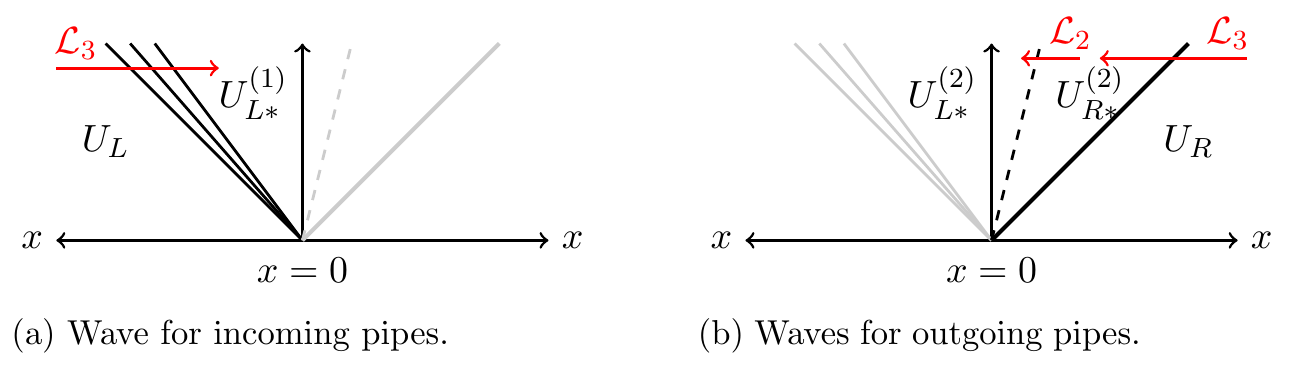}
\parbox{13cm}{
\caption{Connection of the regions $L$, $L\ast$, $R\ast$, and
$R$ with the Lax curve $\LX_3$ for incoming pipes (a) and
the Lax curves $\LX_2\!\circ\!\LX_3$ for outgoing pipes (b).}
\label{fig-lax-curves}
}
\end{figure}
We can now express the coupling conditions for the $\Phi$-solution to the
Riemann problem (\ref{prob_jrm}) in terms of the Lax curves. Remember that
in our network modelling, the $x$-coordinates are chosen in such a way that
pipes are only outgoing from a junction. Consequently, switching from the
standard to the generalised Riemann problem, the sign for the velocity in
incoming pipes has to be changed. This changes the parametrisation of
the $\LX_1$-curve in (\ref{eq_lax1}). A closer inspection of (\ref{eq_rm_u})
shows that $\LX_1$ has to be replaced by $\LX_3$.

Due to the special parametrisation of the pipes and the
restriction to subsonic flow, the contact discontinuity always travels with positive
wave speed and, hence, the state $Y_\ast$ from (\ref{prob_crm}) lies in the region
$L_\ast$, see Fig.~\ref{fig-lax-curves}. We first parameterize all states
$Y^{(i)}_{L\ast}$ using $\LX_3$ for incoming pipes and $\LX_2\circ\LX_3$ for
outgoing pipes, and then apply the function $\Phi$ to them. This yields the set of equations
\begin{equation}
\Phi\left(\left( Y^{(i)}_{L\ast}\right)_{i\in\I_i},
\, \left( Y^{(j)}_{L\ast}\right)_{j\in\I_o}\right) = \bPi\in\R^d
\end{equation}
with
\begin{equation}
\label{eq_lax_star}
Y^{(i)}_{L\ast} = \LX_3(\sigma_i,\bY^{(i)}_0),\;i\in\I_i,\quad\mbox{and}\quad
Y^{(j)}_{L\ast} = \LX_2(\tau_j,\LX_3(\sigma_j,\bY^{(j)}_0)),\;j\in\I_o.
\end{equation}
Let $N_o=\dim(\I_o)$. Then, the degrees of freedom defined
by the Lax curves are $\sigma=(\sigma_1,\ldots,\sigma_N)$ and
$\tau=(\tau_1,\ldots,\tau_{N_o})$. Obviously,
to ensure well-posedness of the generalised Riemann problem at a junction, one
coupling condition has to be provided for incoming pipes, whereas two conditions are
necessary for each of the outgoing pipes. The overall dimension of the parameter space
is $d=N+N_o$.

Given constant states $\bY^{(i)}_0\in D^+,i\!=\!1,\ldots,N_o$, and
$\bY^{(j)}_0\in D^-,j\!=\!N_o+1,\ldots,N$, mass flux, enthalpy and entropy for the $L_\ast$-region can be extracted from formula (\ref{eq_lax_star}):
\begin{equation}
\label{eq_qhs}
\begin{array}{rlll}
f_i(\sigma_i,\tau_i) &=& f_i(\LX_2(\tau_i,\LX_3(\sigma_i,\bY^{(i)}_0))),
                       & i=1,\ldots,N_o, \\
f_j(\sigma_j) &=& f_j(\LX_3(\sigma_j,\bY^{(j)}_0)),
                       & j=N_o+1,\ldots,N,\\
\end{array}
\end{equation}
with $f=q,h,s$. In what follows, we will consider the following coupling
conditions taken from $(M)$, $(H)$, $(S_o)$, and $(S_i)$:
\begin{equation}
\label{eq_coupling}
0 = \Phi(\sigma,\tau) =
\begin{pmatrix}
  \sum_{i=1,\ldots,N_o}\|\nu_i\|\,q_i(\sigma_i,\tau_i)+
  \sum_{j=N_o+1,\ldots,N}\|\nu_j\|\,q_j(\sigma_j)\\[2mm]
  h_{N_o+1}(\sigma_{N_o+1})-h_1(\sigma_{1},\tau_1)\\
  \vdots\\
  h_{N_o+1}(\sigma_{N_o+1})-h_{N_o}(\sigma_{N_o},\tau_{N_o})\\[1mm]
  h_{N_o+1}(\sigma_{N_o+1})-h_{N_o+2}(\sigma_{N_o+2})\\
  \vdots\\
  h_{N_o+1}(\sigma_{N_o+1})-h_{N}(\sigma_{N})\\[1mm]
  s_1(\sigma_1,\tau_1)-s^\ast(\sigma_{N_o+1},\ldots,\sigma_N)\\
  \vdots\\
  s_{N_o}(\sigma_{N_o},\tau_{N_o})-s^\ast(\sigma_{N_o+1},\ldots,\sigma_N)
\end{pmatrix}
\end{equation}
with $s^\ast$ defined through
\begin{equation}
s^\ast = \frac{1}{\sum_{j=N_o+1,\ldots,N}\,\|\nu_j\|q_j(\sigma_j)}
        \,\sum_{j=N_o+1,\ldots,N}\,\|\nu_j\|(q_js_j)(\sigma_j).
\end{equation}
The regularity of the Lax curves ensures the property
$\Phi\in C^1(\R^N\times\R^{N_o},\R^d)$. It remains to show that
(\ref{eq_coupling}) has a unique solution. Then, Newton's method is applied to
determine the solution vector $(\sigma^\ast,\tau^\ast)$, which finally
gives the desired state $Y_\ast$ from
\begin{equation}
Y^{(i)}_\ast = \LX_3(\sigma_i^\ast,\bY^{(i)}_0),\;i\in\I_i,\quad\mbox{and}\quad
Y^{(j)}_\ast = \LX_2(\tau_j^\ast,\LX_3(\sigma_j^\ast,\bY^{(j)}_0)),\;j\in\I_o.
\end{equation}
We note that due to the special choice in (\ref{eq_coupling}) energy and
entropy are conserved at the junction, i.e., $(E)$ and $(S')$ are fulfilled
with $Y_\ast$.

In the case $N=2$ and parallel pipes with the same surface section, the solution
of the generalised Riemann problem coincides with the solution of the standard
Riemann problem for the polytropic Euler equations. We have
\begin{lemma}
Let $N\!=\!2$, $\nu_1\!=\!-\nu_2\ne 0$, and assume constant
initial data $(\bar{\rho}_1,\bar{q}_1,\bar{E}_1)\in D^+$ and
$(\bar{\rho}_2,\bar{q}_2,\bar{E}_2)\in D^-$.
Let $U(x,t)$ be the solution to the standard Riemann problem for
(\ref{euler_eqs}) with initial data
\begin{equation}
\label{eq_crm_data}
U(x,0)=(\rho,q,E)(x,0) =
\begin{cases}
  (\bar{\rho}_1, \bar{q}_1,\bar{E}_1) & \mbox{for }x>0,\\[1mm]
  (\bar{\rho}_2,-\bar{q}_2,\bar{E}_2) & \mbox{for }x<0.
\end{cases}
\end{equation}
Then the functions
\begin{equation}
\begin{array}{rlll}
Y^{(1)}(x,t)=(\rho_1,q_1,E_1)(x,t) &=& (\rho,q,E)(x,t)&\mbox{ if }x>0,\\
Y^{(2)}(x,t)=(\rho_2,q_2,E_2)(x,t) &=& (\rho,-q,E)(-x,t)&\mbox{ if }x<0
\end{array}
\end{equation}
are $\Phi$-solutions in the sense of Def.~\ref{def_phi_sol} that satisfy
the coupling conditions (\ref{eq_coupling}). And vice versa, if
$Y^{(i)}(x,t),i\!=\!1,2,$ are such solutions, then $U(x,t)$ is the solution of
the standard Riemann problem with initial data (\ref{eq_crm_data}).
\end{lemma}
\noindent {\it Proof}: Observe that the assertion holds true if the following
equivalence is satisfied: $\Phi(Y^{(1)}_{L\ast},Y^{(2)}_{L\ast})\!=\!0$ if and only
if $(\rho^{(1)}_{L\ast},q^{(1)}_{L\ast},E^{(1)}_{L\ast})
\!=\!(\rho^{(2)}_{L\ast},-q^{(2)}_{L\ast},E^{(2)}_{L\ast})$. The
coupling conditions simplify to $q_{L\ast}^{(1)}+q_{L\ast}^{(2)}\!=\!0$,
$h_{L\ast}^{(1)}\!=\!h_{L\ast}^{(2)}$, and  $s_{L\ast}^{(1)}\!=\!s_{L\ast}^{(2)}$.
Since the solution is smooth along $x\!=\!0$, density and total energy are uniquely
determined by the values of $h$ and $s$. This gives the desired equality.\qedwhite \\

\noindent For the general case of $N$ connected pipes at one junction, we can
show a local result for the well-posedness of the generalised Riemann problem
(\ref{prob_jrm}) with the coupling function $\Phi$ defined through $(M)$, $(H)$,
$(S_o)$, $(S_i)$ and stated in more detail in (\ref{eq_coupling}). Similar results
can be found in \cite[Theorem 2.7]{ColomboMauri2008} and \cite[Proposition 2.4]{Herty2008}
for other coupling conditions.
\begin{theorem}
\label{th_riemann}
Let $N>N_0>0$ and $\Phi$ defined through $(M)$, $(H)$, $(S_o)$, and $(S_i)$.
Assume constant initial data
$\bU^{(i)}\in D^+,i\!=\!1,\ldots,N_o$, and
$\bU^{(j)}\in D^-,j\!=\!N_o+1,\ldots,N$,
with $\Phi(\bU)\!=\!0$ are given. Then there exist positive constants $\delta$ and $K$
such that for all initial states $\tU\in(\R^+\times\R\times\R^+)^N$ with
$\sum_{i=1,\ldots,N}\|\tU^{(i)})-\bU^{(i)})\|\!<\!\delta$, the Riemann problem (\ref{prob_jrm}) admits
a unique $\Phi$-solution $U(x,t)=\cR^{\Phi}(\tU)$ satisfying $\Phi(U(0^+,t))\!=\!0$ and
\begin{equation}
\label{grm_lip_init}
\|\cR^{\Phi}(\tU)-\cR^{\Phi}(\bU)\|_{{\bf L}^\infty(\Omega)} \le K\,
\sum_{i=1}^N\|\tilde{U}^{(i)}-\bU^{(i)}\|.
\end{equation}
Additionally, if $\nu$ is replaced by $\hat{\nu}$, where
$\sum_{i=1,\ldots,N}\|\nu_i-\hat{\nu}_i\|\!<\!\delta$, and
$\cR^{\Phi}_{\hat{\nu}}(\tU)$ is the corresponding $\Phi$-solution
for the same initial state $\tU$, then
\begin{equation}
\label{grm_lip_pipe}
\|\cR^{\Phi}_\nu(\tU)-\cR^{\Phi}_{\hat{\nu}}(\tU)\|_{{\bf L}^\infty(\Omega)}
\le K\,\sum_{i=1}^N\|\nu_i-\hat{\nu}_i\|
\end{equation}
with $\cR^{\Phi}_\nu(\tU)\!=\!\cR^{\Phi}(\tU)$.
\end{theorem}
\noindent {\it Proof}: We follow the proof of Theorem 2.7 in \cite{ColomboMauri2008}
and show that (\ref{eq_coupling}) has locally a unique
solution. Observe $\Phi(\sigma,\tau)\!=\!0$ for $\sigma_0\!=\!(\bar{p}_1,\ldots,\bar{p}_N)$, and $\tau_0\!=\!0\in\R^{N_o}$, since the initial data satisfy the coupling conditions.
In the spirit of the implicit function theorem, it is
sufficient to study the determinant of the Jacobian $D_{(\sigma,\tau)}\Phi(\sigma_0,\tau_0)$.

Let us first collect a few derivatives. For incoming pipes, we derive from the second equation in (\ref{eq_qhs})
\begin{equation}
q_j'(\bar{p}_j)=\frac{\lambda_3(\bar{u}_j)}{\bar{c}_j^2},\quad
h_j'(\bar{p}_j)=\frac{\lambda_3(\bar{u}_j)}{\bar{c}_j\bar{\rho}_j},\quad
\partial_{\sigma_j}s^\ast(\bar{p})=\frac{\|\nu_j\|\lambda_3(\bar{u}_j)}{\bar{c}_j^2
\sum_{i\in\I_i}\|\nu_i\|\bar{q}_i}\left( \bar{s}_j-\bar{s}^\ast\right)
\end{equation}
with $\bar{c}_j=\sqrt{\gamma\bar{p}_j/\bar{\rho}_j}$ and $j=N_o+1,\ldots,N$. Further,
the first equation in (\ref{eq_qhs}) yields for outgoing pipes
\begin{eqnarray}
\partial_{\sigma_i}q_i(\bar{p}_i,0)=\frac{\lambda_3(\bar{u}_i)}{\bar{c}_i^2},\quad
\partial_{\sigma_i}h_i(\bar{p}_i,0)=\frac{\lambda_3(\bar{u}_i)}{\bar{c}_i\bar{\rho}_i},\quad
\partial_{\sigma_i}s_i(\bar{p}_i,0)=0,\\
\label{eq_with_gamma}
\partial_{\tau_i}q_i(\bar{p}_i,0)=\lambda_2(\bar{u}_i),\quad
\partial_{\tau_i}h_i(\bar{p}_i,0)=-\frac{\bar{c}_i^2}{(\gamma-1)\bar{\rho}_i},\quad
\partial_{\tau_i}s_i(\bar{p}_i,0)=-\frac{\gamma c_v}{\bar{\rho}_i}
\end{eqnarray}
for $i=1,\ldots,N_o$. This yields the following matrix for the Jacobian $D_{(\sigma,\tau)}\Phi(\sigma_0,\tau_0)$:
\begin{equation}
\label{cp_jac}
\left(
\begin{array}{ccc|cccc|ccc}
        \hat{q}_{\sigma_1} &
            \!\!\cdots\!\! &
    \hat{q}_{\sigma_{N_o}} &
  \hat{q}_{\sigma_{N_o+1}} &
  \hat{q}_{\sigma_{N_o+2}} &
            \!\!\cdots\!\! &
      \hat{q}_{\sigma_{N}} &
          \hat{q}_{\tau_1} &
            \!\!\cdots\!\! &
        \hat{q}_{\tau_{N_o}}\\[1mm] \hline
       -h_{\sigma_1} & & &
  h_{\sigma_{N_o+1}} & & & &
         -h_{\tau_1} & &\\
 & \!\!\ddots\!\! & & \vdots & & & & & \!\!\ddots\!\! &\\
                     & &
   -h_{\sigma_{N_o}} &
  h_{\sigma_{N_o+1}} & & & & & &
       -h_{\tau_{N_o}}\\[1mm] \hline
                     & & &
  h_{\sigma_{N_o+1}} &
 -h_{\sigma_{N_o+2}} & & & & &\\
  & & & \vdots & & \!\!\ddots\!\! & & \\
                     & & &
  h_{\sigma_{N_o+1}} & & &
     -h_{\sigma_{N}} & & &\\[1mm] \hline
                     & & &
  -s^\ast_{\sigma_{N_o+1}} &
  -s^\ast_{\sigma_{N_o+2}} &
            \!\!\cdots\!\! &
  -s^\ast_{\sigma_{N}} &
            s_{\tau_1} & &\\
 & & & \vdots & & & \vdots & & \!\!\ddots\!\! &\\
                      & & &
  -s^\ast_{\sigma_{N_o+1}} &
  -s^\ast_{\sigma_{N_o+2}} &
            \!\!\cdots\!\! &
      -s^\ast_{\sigma_{N}} & & &
            s_{\tau_{N_o}}
\end{array}
\right)
\end{equation}
Here, we have used the short notations $f_{\mu_i}=\partial_{\mu_i}f_{i}$,
$\hat{q}_{\mu_i}=\|\nu_i\|\partial_{\mu_i}q_{i}$
for $f\!=\!h,s$, and $\mu=\sigma,\tau$, and
$s^\ast_{\sigma_{i}}=\partial_{\sigma_{i}}s^\ast$. Observe that none of
the derivatives can vanish, except $s^\ast_{\sigma_{i}}$. We find
\begin{equation}
\hat{q}_{\sigma_i}>0,\;\hat{q}_{\tau_j}>0,\;h_{\sigma_i}>0,\;h_{\tau_j}<0,\;s_{\tau_i}<0
\quad\mbox{for }i=1,\ldots,N,\;j=1,\ldots,N_o.
\end{equation}
Without loss of generality, we choose the numbering of the incoming pipes in such
a way that $\bar{s}_{N_o+1}=\max_{i\in\I_i}\bar{s}_i$. Then
$\bar{s}_{N_o+1}-\bar{s}^\ast\ge 0$, and since $\bar{q}_i\!<\!0$ for
$i\in\I_i$, it follows that $s^\ast_{\sigma_{N_o+1}}\le 0$. From the special
structure of the matrix (\ref{cp_jac}), we deduce that the Jacobian is regular if
and only if all $3\times 3-$matrices
\begin{equation}
D_i=
\begin{pmatrix}
  \hat{q}_{\sigma_i} & \hat{q}_{\sigma_{N_o+1}} & \hat{q}_{\tau_i}\\[1mm]
  -h_{\sigma_i} & h_{\sigma_{N_o+1}} & -h_{\tau_i}\\[1mm]
  0 & -s^\ast_{\sigma_{N_o+1}} & s_{\tau_i}
\end{pmatrix}
\quad\quad\mbox{for }i=1,\ldots,N_o,
\end{equation}
are regular. Taking into account the signs of all derivatives, we have
\begin{equation}
\det \left( D_i \right) = \hat{q}_{\sigma_i}
(h_{\sigma_{N_o+1}}s_{\tau_i}-h_{\tau_i}s^\ast_{\sigma_{N_o+1}})+
h_{\sigma_i} (\hat{q}_{\sigma_{N_o+1}}s_{\tau_i}+
\hat{q}_{\tau_i}s^\ast_{\sigma_{N_o+1}}) < 0.
\end{equation}
Therefore, $\det(D_{(\sigma,\tau)}\Phi(\sigma_0,\tau_0))\ne 0$ and by the
implicit function theorem, there exist a $\delta\!>\!0$, a neighbourhood
$\cU(v_0)$  of $v_0=(\sigma_0,\tau_0)$, and a function
$\varphi:B(\bar{U},\delta)\rightarrow \cU(v_0)$ such that $\varphi(\bar{U})=v_0$
and $\Phi(v;U)=0$ if and only if $v=\varphi(U)$ for all $U\in B(\bar{U},\delta)$.
The solution $U(x,t)$ can then be identified by the restriction to $x\in\R^+$ of
the solution to the standard Riemann problem (\ref{prob_crm}) with $\bY_0\!=\!\tU$
and
\begin{equation}
Y_\ast^{(i)}=\LX_3(\varphi(\tU)_i,\tU),\;i\in\I_i,\mbox{ and }
Y_\ast^{(j)}=\LX_2(\varphi(\tU)_{j+N},\LX_3(\varphi(\tU)_j,\tU),\;j\in\I_o.
\end{equation}
The Lipschitz estimate (\ref{grm_lip_init}) follows from the $C^1$-regularity of $\Phi$.
Since $\Phi$ depends smoothly on $\|\nu_i\|$, the same arguments as above can be used to
show (\ref{grm_lip_pipe}).
\qedwhite

\begin{remark}(energy and entropy conservation) We would like to remember that the
coupling conditions ensure conservation of energy and entropy at the junction,
\begin{equation}
\sum_{i=1}^{N}\,\|\nu_i\|\,(u_i(E_i+p_i))(0^+,t) =
\sum_{i=1}^{N}\,\|\nu_i\|\,(q_is_i)(0^+,t) = 0.
\end{equation}
It is therefore not necessary to assume that the perturbed initial state
$\tU$ is strictly entropic, i.e., satisfies the strict entropy inequality
in $(S)$ as used in \cite{ColomboMauri2008,Herty2008}.
\end{remark}

\begin{remark}
Theorem \ref{th_riemann} remains valid even if the adiabatic
exponent $\gamma$ varies over the set of pipes. In this case,
$\bar{c}_i=\sqrt{\gamma_i\bar{p}_i/\bar{\rho}_i}$ and $\gamma$
has to be replaced by an individual $\gamma_i>1$ in (\ref{eq_with_gamma}),
which does not influence the sign arguments
used in the proof.
\end{remark}

\section{The Cauchy problem at the junction}
\label{sec:cauchy_prob}
In this section, we define a weak entropic solution for the
general Cauchy problem with source terms at junctions, using
the above stated coupling conditions. Further, two main results
are formulated: the well-posedness for the homogeneous as well as the
inhomogeneous case under the well known assumption
that the total variation of the initial data is sufficiently small.
Both theorems can be seen in line with
Theorem~3.2. from {\sc Colombo and Mauri} \cite{ColomboMauri2008}
and Theorem~2.3. from {\sc Colombo, Guerra, Herty, and
Schleper} \cite{ColomboGuerraHertySchleper2009}. The key point
is the well-posedness of the Riemann problem stated in
Theorem~\ref{th_riemann} above,
which provides the basis for the proofs.

We first introduce a few notations.
\begin{definition}
Let
\begin{equation}
\begin{array}{rcll}
\|Y\| &=& \ds\sum_{i=1}^N \left\| Y^{(i)}\right\| & \mbox{for } Y\in\Omega,\\[5mm]
\|Y\|_{\bfL^1} &=& \ds\int_{R^+}\|Y(x)\|\,dx & \mbox{for } Y\in\bfL^1(R^+;\Omega)\\[5mm]
TV(Y) &=& \ds\sum_{i=1}^N TV(Y^{(i)}) & \mbox{for } Y\in\bfBV(R^+;\Omega).
\end{array}
\end{equation}
For a constant state $\bY$ and a positive $\delta\in[0,\bdelta]$, we set
\begin{equation}
D_\delta(\bY) = \{ Y \in \bY + \bfL^1(R^+;\Omega):\;TV(Y)\le \delta\}.
\end{equation}
\end{definition}
Let $G$ denote the vector of the right-hand side functions in
(\ref{euler_eqs}) for all pipes and be defined through
\begin{equation}
(G(t,Y))(x) = \left( G(x,t,Y^{(1)}),\ldots,G(x,t,Y^{(N)})\right).
\end{equation}
For the map $G:[0,T]\times D_{\bdelta}(\bY)\rightarrow \bfL^1(R^+;\Omega)$,
we assume that there exist positive constants $L_1$ and $L_2$ such that
for all $t,s\in [0,T]$ the following inequalities are satisfied:
\begin{equation}
\begin{array}{rlll}
\|G(t,Y_1)-G(s,Y_2)\|_{\bfL^1} &\!\le\!& L_1 \left(\|Y_1-Y_2\|_{\bfL^1}
+ |t-s| \right) & \mbox{for all }Y_1,Y_2\in D_{\bdelta}(\bY),\\[2mm]
TV(G(t,Y)) &\!\le\!& L_2 & \mbox{for all }Y\in D_{\bdelta}(\bY).
\end{array}
\end{equation}
This is the usual assumption on $G$, which also covers non-local terms
\cite{ColomboGuerra2007,ColomboGuerra2008} as well as real applications
\cite{ColomboGuerraHertySchleper2009}.

Next we define the Cauchy problem at junctions, which corresponds to our
special set of coupling conditions.
\begin{definition}
\label{def_sol_cauchy}
Let $N>N_0>0$ and $\Phi$ defined through $(M)$, $(H)$, $(S_o)$, and $(S_i)$.
A weak solution on $[0,T]$ to the Cauchy problem
\begin{equation}
\label{prob_cauchy}
\begin{array}{rcll}
\partial_t U^{(i)} + \partial_x F(U^{(i)}) &=& G(x,t,U^{(i)}), &
(x,t)\in\R^+\times\R^+,i=1,\ldots,N,\\[1mm]
\Phi(U(0^+,t)) &=& 0,& t\in\R^+,\\[1mm]
U(x,0) &=& U_0(x), & x\in \R^+,\;U_0\in\bU+\bfL^1(\R^+;\Omega),
\end{array}
\end{equation}
is a map $U\in\bfC^0([0,T];\bU+\bfL^1(\R^+;\Omega))$ that corresponds
to $\bfBV(\R^+;\Omega)$ for all $t\in [0,T]$ and
satisfies the initial condition, $U(x,0)\!=\!U_0(x)$,
and the condition at the junction, $\Phi(U(0^+,t))=0$,
for a.e. $t>0$. Further,
for all $\varphi\in\bfC^{\infty}_c(\R^+\times (0,T);\R)$ it holds
\begin{equation}
\ds\sum_{i=1}^N\left(\int_0^T \int_{\R^+} \left( \rho_i\partial_t\varphi +
q_i\partial_x\varphi + G_1(x,t,U^{(i)})\varphi\right)\,dx\,dt
\right) \,\|\nu_i\| = 0
\end{equation}
and
\begin{equation}
\begin{array}{rcl}
\ds\int_0^T \int_{\R^+} \left( q_i\partial_t\varphi +
P_i\partial_x\varphi + G_2(x,t,U^{(i)})\varphi\right)\,dx\,dt &=&
\ds\int_0^T P_i(0^+,t)\varphi(0,t)\,dt,\\[5mm]
\ds\int_0^T \int_{\R^+} \left( E_i\partial_t\varphi +
q_ih_i\partial_x\varphi + G_3(x,t,U^{(i)})\varphi\right)\,dx\,dt &=&
\ds\int_0^T q_i(0^+,t)h^\ast(t)\varphi(0,t)\,dt.
\end{array}
\end{equation}
for all $i=1,\ldots,N$ with $P_i=\rho_iu_i^2+p_i$ and a suitable
$h^\ast(t)\in L^1([0,T];\R^+)$.

The weak solution is entropic if for all non-negative $\varphi\in\bfC^{\infty}_c(\Ro^+\times(0,T);\R^+)$
and $i=1,\ldots,N$
\begin{equation}
\int_0^T \int_{\R^+} \left( \rho_is_i\partial_t\varphi +
q_is_i\partial_x\varphi +
\partial_U(\rho_is_i)G(x,t,U^{(i)})\varphi\right)\,dx\,dt \ge 0.
\end{equation}
\end{definition}
We note that multiplying the energy equation with $\|\nu_i\|$ and summing
up over all pipes gives the energy balance equation
\begin{equation}
\ds\sum_{i=1}^N\left(\int_0^T \int_{\R^+} \left( E_i\partial_t\varphi +
q_ih_i\partial_x\varphi + G_3(x,t,U^{(i)})\varphi\right)\,dx\,dt
\right) \,\|\nu_i\| = 0,
\end{equation}
which means energy conservation in the case $G_3=0$.

A solution to the Cauchy problem can be constructed by means of the wave front
tracking method. In the book of {\sc Bressan} \cite{Bressan2000} all
necessary steps can be found.

Let us first consider the homogeneous case. We have the following
\begin{theorem}
Let $G\!=\!0$, $N>N_0>0$ and $\Phi$ defined through $(M)$, $(H)$, $(S_o)$, and $(S_i)$.
Assume constant initial data $\bU^{(i)}\in D^+,i\!=\!1,\ldots,N_o$, and
$\bU^{(j)}\in D^-,j\!=\!N_o+1,\ldots,N$, with $\Phi(\bU)\!=\!0$ are given.
Then there exist positive constants $\delta$, $K$, and a uniformly Lipschitz semigroup
$S:\R^+\times D\rightarrow D$ such that:
\begin{itemize}
\item[(1)] $\overline{D_\delta}(\bU)\subseteq D$.
\item[(2)] $S_0=Id$ and $S_sS_t=S_{s+t}$.
\item[(3)] For all $U\in D$, the map $t\rightarrow S_t(U)$ is a weak
entropic solution to the Cauchy problem (\ref{prob_cauchy}) in the sense
of Definition \ref{def_sol_cauchy}.
\item[(4)] For $\hU,\tU\in D$ and $s,t\ge 0$
\[
\|S_t(\hU)-S_s(\tU)\|_{\bfL^1(\R^+;\Omega)}\le
K\,(\|\hU-\tU\|_{\bfL^1(\R^+;\Omega)}+|s-t|).
\]
\item[(5)] If $U\in D$ is piecewise constant and $t>0$ sufficiently small,
then $S_t(U)$ coincides with the juxtaposition of the solutions to Riemann
problems centered at the points of jumps or at the junction.
\end{itemize}
\end{theorem}
\noindent {\it Proof}: The properties are a direct consequence of a natural
extension of the standard Riemann semigroup theory
\cite[Section 8.3]{Bressan2000} to junctions. All arguments can be copied
from the proof of Theorem 3.2. in \cite{ColomboMauri2008}. \qedwhite\\

\noindent For non-vanishing sources $G$, we get the following result for the
well-posedness of the Cauchy problem:
\begin{theorem}
Let $N>N_0>0$ and $\Phi$ defined through $(M)$, $(H)$, $(S_o)$, and $(S_i)$.
Assume constant initial data $\bU^{(i)}\in D^+,i\!=\!1,\ldots,N_o$, and
$\bU^{(j)}\in D^-,j\!=\!N_o+1,\ldots,N$, with $\Phi(\bU)\!=\!0$ are given.
Then there exist positive constants $\delta$, $\delta'$, $K$, domains
$D_t$ for $t\in [0,T]$, and a map $\cE(s,t_0):D_{t_0}\rightarrow D_\delta$
with $t_0\in[0,T]$ and $s\in [0,T-t_0]$ such that
\begin{itemize}
\item[(1)] $D_{\delta'}(\bU)\subseteq D_t\subseteq D_{\delta}(\bU)$
for all $t\in[0,T]$.
\item[(2)] $\cE(0,t_0)U=U$ for all $t_0\in [0,T],\;U\in D_t$.
\item[(3)] $\cE(s,t_0)D_{t_0}\subset D_{t_0+s}$ for all
$t_0\in [0,T],\;s\in [0,T-t_0]$.
\item[(4)] For all $t_0\in [0,T]$, $s_1,s_2\ge 0$
with $s_1+s_2\in [0,T-t_0]$
\[ \cE(s_2,t_0+s_1)\circ\cE(s_1,t_0)=\cE(s_1+s_2,t_0). \]
\item[(5)] For all $U_0\in D_{t_0}$, the map $t\rightarrow\cE(t,t_0)U_0$
is the entropic solution to the Cauchy problem (\ref{prob_cauchy})
in the sense of Definition \ref{def_sol_cauchy}.
\item[(6)] For all $t_0\in [0,T]$ and $U_0\in D_{t_0}$
\[ \lim_{t\rightarrow 0}\frac{1}{t}
\|U(t)-(S_t(U_0)+tG(t_0,U_0))\|_{\bfL^1}=0,\]
where $U(t)=\cE(t,t_0)U_0$ and $S_t$ denotes the semigroup generated from
(\ref{prob_cauchy}) with $G=0$.
\item[(7)] For all $t_0\in [0,T],\;s\in [0,T-t_0]$ and
$U,\tU\in D_{t_0}$
\[ \| \cE(s,t_0)U-\cE(s,t_0)\tU\|_{\bfL^1} \le K \|U-\tU\|_{\bfL^1}.\]
\end{itemize}
\end{theorem}
\noindent {\it Proof}: The proof can be achieved by following the standard
line developed in
\cite{ColomboGuerraHertySchleper2009} for $2\times 2$ hyperbolic systems.
We set $\Pi=0$ and use a modified version of the Glimm type and Bressan-Liu-Yang
functionals,
which are obtained by an extension to the present case of a $3\times 3$
Euler system by means of
the techniques presented in \cite{Bressan2000,DonadelloMarson2007}. This is
straightforward and bears no difficulties. \qedwhite

\section{Summary}
\label{sec:summary}
We have proposed a novel set of physically sound coupling conditions
at a junction of pipes with possibly different cross-sectional areas
for the $3\times 3$ one-dimensional system of homogeneous Euler equations.
In the subsonic flow regime, these conditions ensure mass, energy and entropy
conservation at the junction. The new approach is
applicable for general situations with at least one incoming and one outgoing
pipe. Previously used pressure-based coupling conditions that can produce non-physical
solutions are replaced by physically sound entropy-preserving conditions. The equality
of the entropy at the junction-pipe interface for pipes with outgoing flow is
enforced and the corresponding coupling constant is identified as the convex
combination of all entropies that belong to the pipes with incoming flow.
The existence and uniqueness of solutions to generalised Riemann problems
at a junction in the neighbourhood of constant stationary states are proven.
Following standard proof techniques, this yields the well-posedness of the
homogeneous and inhomogeneous Cauchy problems for initial data with
sufficiently small total variation.

\section{Acknowledgement}
This work was supported by the
German Research Foundation within the collaborative research center
TRR154 ``Mathematical Modeling, Simulation and Optimization Using
the Example of Gas Networks'' (DFG-SFB TRR154/1-2014, TP B01).

\bibliographystyle{plain}
\bibliography{bibeuler}

\end{document}